\documentclass[reqno]{amsart}
\usepackage{latexsym,amssymb}
\usepackage[bottom]{footmisc}
\usepackage{graphicx}
\usepackage{xspace}
\usepackage[bookmarksnumbered,colorlinks]{hyperref}
\usepackage{graphics}

\newcommand{\bt}{\begin{teo}}                              
\newcommand{\et}{\end{teo}}                                
\newcommand{\bco}{\begin{cor}}                              
\newcommand{\eco}{\end{cor}}                                
\newcommand{\bd}{\begin{defn}}                           
\newcommand{\ed}{\end{defn}}                             
\newcommand{\bl}{\begin{lem}}                            
\newcommand{\el}{\end{lem}}                              
\newcommand{\bpr}{\begin{prop}}               
\newcommand{\epr}{\end{prop}}                 
\newcommand{\bere}{\begin{remark}}            
\newcommand{\ere}{\end{remark}}               
\newcommand{\beq}{\begin{equation}}
\newcommand{\eeq}{\end{equation}}
\def\bal#1\eal{\begin{align}#1\end{align}}                      
\def\baln#1\ealn{\begin{align*}#1\end{align*}}         
\def\bml#1\eml{\begin{multline}#1\end{multline}}       
\def\bmln#1\emln{\begin{multline*}#1\end{multline*}}  
\def\bga#1\ega{\begin{gather}#1\end{gather}}
\def\bgan#1\egan{\begin{gather*}#1\end{gather*}}
\newcommand{\de}{\mathrm{d}}

\newcommand{\R}{\ensuremath{\mathbb{R}}\xspace}     
                      
\newtheorem{teo}{Theorem}[section]

\newtheorem{lem}[teo]{Lemma}
\newtheorem{defn}[teo]{Definition}
\newtheorem{cor}[teo]{Corollary}
\newtheorem{prop}[teo]{Proposition}
\theoremstyle{remark}
\newtheorem{remark}[teo]{Remark}
\theoremstyle{remark}
\newtheorem{problem}[teo]{Problem}

\title[Convexity of a hypersurface in a semi-Riemannian manifold]{Infinitesimal and local convexity of a hypersurface in a semi-Riemannian manifold}

\author[E. Caponio]{Erasmo Caponio}
\address{Dipartimento di Meccanica, Matematica e Managment\hfill\break\indent
Politecnico di Bari \hfill\break\indent Via Orabona 4,
70125, Bari, Italy}
\email{caponio@poliba.it}



\date{}

\begin{document}

\begin{abstract}
Given a Riemannian manifold $(M,g)$ and an embedded hypersurface $H$ in $M$, a result by R. L. Bishop states that  infinitesimal convexity on a neighborhood of a point in $H$ implies local convexity. Such result was extended very recently to Finsler manifolds by the author et al. \cite{BaCaGS11}. We show in this note that the techniques in \cite{BaCaGS11}, unlike the ones in Bishop's paper, can be used to prove the same result when $(M,g)$ is  semi-Riemannian. We make some remarks for the case when only timelike, null or spacelike geodesics are involved. The notion of geometric convexity is also reviewed and some applications to geodesic connectedness of an open subset of a Lorentzian manifold are given.
\end{abstract}
\maketitle

\section{Introduction}
Let $\Omega$ be an open of subset of $\R^2$ and $f\colon \Omega\to \R$ be a differentiable function in $\Omega$. We say that $f$ is locally convex at a point $(x_0,y_0)\in \Omega$ if there exists a neighborhood $U\subset \Omega$ of  $(x_0,y_0)$ such that the graph of the restriction of $f$ to $U$ is nowhere below the tangent plane at $\big(x_0,y_0, f(x_0,y_0)\big)$ or, equivalently, all the straight lines through the point $\big(x_0,y_0,f(x_0,y_0)\big)$ are  nowhere above the graph of $f_{|U}$. If  $f$ is twice differentiable then the above condition is satisfied if  the Hessian of $f$ is positive semidefinite in a neighborhood of $(x_0,y_0)$, while the positive semi-definiteness of $\de^2 f$ at the single point $(x_0,y_0)$ is necessary but not sufficient, as the function $f(p)=f(x,y)=-(x^4+y^4)$ at $(0,0)$ shows.

Let $H$ be the hypersurface in $\R^3$ which is the graph of the function $f$; the  following quadratic form associated to the Hessian of $f$ at a point $(x,y)$, i.e. 
\[\Pi_{(x,y)}\big((u,v),(u,v)\big)= \frac{f_{xx}(x,y)u^2+2f_{xy}(x,y)uv+f_{yy}(x,y)v^2}{\sqrt{1+f^2_x(x,y)+f^2_y(x,y)}},\] 
is the {\em second fundamental form} of $H$ at $\big(x,y,f(x,y)\big)$. Thus, the convexity of a twice differentiable function $f$ in a neighborhood of $(x_0,y_0)\in\Omega$ is equivalent to the fact that the $\Pi$ is positive semidefinite in the same neighborhood.

Now assume that $H$ is a smooth embedded  hypersurface in a Riemannian manifold $(M,g)$. The natural  generalizations of the above notions are the following ones:
\begin{itemize}
 \item[$\triangledown$]the first one becomes the requirement that there exists a neighborhood $\tilde U$ in $M$ of $p_0\in H$ such that the intersections of the images of the geodesics through $p_0$ with velocity vector at $p_0$ tangent to $H$ with $\tilde U$ are contained on the closure of one of the two connected component $\tilde U\setminus H$ or, equivalently, there exists a neighborhood $O$ of $0\in T_{p_0} H$ and a neighborhood $\tilde U$ of $p_0$ such that the $\exp_{p_0}(O)$ is contained $H\cup C$, where $C$ is one of the two connected component of $\tilde U\setminus H$; this condition is called {\em local convexity} at $p_0$;
\item[$\Diamond$]\label{infconv}the second one becomes the requirement that there exists an open neighborhood  $U\subset H$ of $p_0$ such that the second fundamental form  of $U$, with respect to a smooth unit normal vector field $n$ on $U$, 
is positive semidefinite\footnote{According to the  case when $H$ is  the graph of a function and  $n$ is the smooth unit vector which lies on the same side of the graph as the canonic unit vector defining the axis of the values of the function, from \cite{Bishop75}, we know that the geodesics issuing from $p_0$ and tangent to $H$  are locally contained in the closure of the component of a tubular neighborhood of $H$  individuated by $-n$. As the sign of $\Pi$ changes if one changes $n$ with $-n$, it is clear that what is really important in this definition is the fact that $\Pi$ is semidefinite in $U$, either positive or negative. Clearly in the case when $\Pi$ is negative semidefinite the geodesics tangent at $p_0$ to $H$ are locally contained in the component individuated by $n$ itself. \label{nota}}  for each $p\in U$, i.e.  
\beq\label{infconvineq}
\Pi_p(v,v)=g_p(n,\nabla_VV)\geq 0,
\eeq 
where $\nabla$ is the Levi-Civita connection of $(M,g)$ and $V$ is a vector field on $M$ extending $v\in T_p H$; this condition is called {\em infinitesimal convexity} in $U$.
\end{itemize}
\bere
The global versions of these definitions require that the embedded hypersurface $H$ is {\em orientable}, i.e. $H$ must admits a smooth unit normal vector field $n$.  In such a case, $H$ is said infinitesimally convex if $\Pi(v,v)\geq 0$, for each $v\in TH$, where $\Pi$ is the second fundamental form of $H$ with respect to $n$.
\ere
As in the case of the graph  of a function, infinitesimal convexity at a neighborhood of  $p_0\in H$ implies local convexity. Surprisingly enough, this property was proved only relatively recently for a hypersurface in a Riemannian manifold having constant sectional curvature in \cite{CarWar70} and in the general case in \cite{Bishop75}. Besides, the proof in \cite{Bishop75} is rather involved and it seems to depend drastically on the  Riemannian setting. In Section~\ref{bishop}, we will describe which are the obstacles one encounters trying to extend Bishop's proof when the ambient manifold is a Finsler or a semi-Riemannian one. The implication in the Finsler setting was obtained recently in \cite{BaCaGS11}. In Section~\ref{main} we will show that the techniques in \cite{BaCaGS11} also give the implication in a  semi-Riemannian manifold and in Section~\ref{timespace} we will consider the same problem in relation to the causal character of the geodesics and to another notion of convexity called {\em 
geometric 
convexity}. Finally, we will see some applications to geodesic connectedness  of an open subset of a Lorentzian manifold by means of 
causal geodesics. 
\section{A review of Bishop's proof}\label{bishop}
The proof in \cite{Bishop75} that infinitesimal convexity in a neighborhood $U$ of a point $p_0\in H$ implies local convexity at $p_0$  is based on a reduction of the general problem to the case of a Riemannian manifold $(M, g)$ of dimension two. Such a reduction is quite natural and it is based on the following idea. For any geodesic $\gamma$ starting from $p_0$ with initial velocity tangent to $H\cap U$ at $p_0$, consider the surface $S(\gamma)$ ruled  by the geodesics orthogonal to $H\cap U$ and passing through $\gamma$. Two such surfaces $S(\gamma)$, $S(\bar \gamma)$, $\gamma\neq \bar \gamma$, intersect along the geodesic $\sigma$ which is orthogonal to $H$ at $p_0$ 
(see Fig. 1).
\begin{figure}[h]
\includegraphics[scale=1.4]{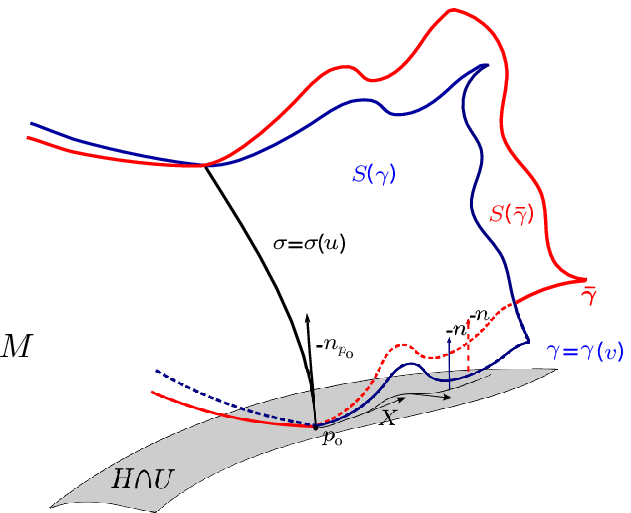}
%
%
\caption{}
\end{figure}
Then, for each $\gamma$, define the unit vector field $E\in TS(\gamma)$ which is the geodesic field (with respect to the Riemannian metric induced on $S(\gamma)$  by $g$)
such that $E_p$ is orthogonal to $\sigma$ for all $p\in\sigma$. The surface $S(\gamma)$ can be parametrized by the coordinates  $u$ and $v$ which are, respectively, the affine parameters of $\sigma$ and $\gamma$ (thus $\gamma$ is represented on $S(\gamma)$ by the equation $u=0$).  The unit vector field $X$ in $TS(\gamma)$ which is tangent to $H\cap S(\gamma)$ can be written as $X_{(u,v)}=c(u,v)E_{(u,v)} + s(u,v)F_{(u,v)}$, where $F$ is the smooth unit vector field  in $TS(\gamma)$ orthogonal to $E$ (oriented in such way that $F_{p_0}=-n_{p_0}$). Using the infinitesimal convexity assumption, it is possible to obtain a differential inequality, satisfied by  the function $s$, implying  that $s$  is non positive on $H\cap S(\gamma)$. Then local convexity at $p$ follows from the observation that $X(u)=cE(u)+sF(u)=sF(u)\geq 0 $. 

The delicate point in this argument is that $E$ could not be a well defined vector field in $\tilde U\setminus \{\sigma\}$, for any convex  neighborhood $\tilde U$ of $p_0$ in $M$. 

By {\em convex neighborhood} of a point $p$ in a Riemannian or semi-Riemannian manifold,  here and hereafter, we mean an open  neighborhood $\tilde U$ which has the property that any two point in $\tilde U$ are connected by a unique geodesic whose support is contained in $\tilde U$ (it is well known than any point on a Riemannian or semi-Riemannian manifold has a convex neighborhood, see for example \cite[Ch. 5, Prop. 7]{One83}).   

Bishop rules out this pathological case by using an estimate on the distance in $M$ of the focal points of $\sigma$ in $S(\gamma)$. Such estimate is based on the one hand on the fact that,   the Gauss curvatures of the surfaces $S(\gamma)$ are uniformly bounded  from above (in a neighborhood of $p_0$) by the sectional curvature of the ambient manifold $M$ (provided that the metric is at least $C^4$), on the other one on the fact that the geodesics in $S(\gamma)$ having $E$ as their velocity vector fields, leave $\sigma$ orthogonally and then the the distance from $\sigma$ (in $M$) of their points  can be evenly controlled. 

Such type of arguments have no straightforward extension  when the ambient manifold $(M,g)$ is semi-Riemannian. 
Indeed, as shown in \cite{Nomizu83}, if the sectional curvature is bounded for all timelike or spacelike planes at a point  $p$ then it is constant at $p$. Moreover the Morse index of spacelike geodesics is always $+\infty$. Thus,  Rauch comparison theorems (see \cite{Berger62} for the case of geodesics starting orthogonally to a given geodesic)   do not generally hold in semi-Riemannian geometry. Some results are available in literature, but  they concern either causal geodesics (cf. \cite[\S 11.2]{BeErEa96} for geodesics connecting two given points and \cite{Kim90,EhrKim91} for more general boundary conditions) or Jacobi fields  along a given   single geodesic of any causal character satisfying two different types of boundary conditions \cite{JavPic09}.

In \cite{BaCaGS11},  Bishop's result was extended to a hypersurface in a Finsler manifold $(M,F)$
by using a different  approach from that in \cite{Bishop75}. In the Finslerian setting, the problem of directly extending the proof in \cite{Bishop75} is not related to the lack of satisfactory  comparison results (cf. \cite[\S 9.1]{BaChSh00}) but rather
on the fact that the differential inequality satisfied by the function $s$ above (that comes from the evaluation of the second fundamental form of the curve $S(\gamma)\cap H$, in the two dimensional manifold $S(\gamma)$, which assumes by construction the same values of that of the hypersurface $H$) is obtained by using, in an essential way, the metric compatibility of the Levi-Civita connection. For a Finsler manifold, the role of the Levi-Civita connection can be played  by the Chern connection which is metric compatible (and torsion free) if and only the Finsler structure is Riemannian (see, e.g. \cite[Exercise 2.4.2]{BaChSh00}).\footnote{In some cases, the Chern connection which is a linear connection on the vector bundle $\pi^*TM$ over $TM\setminus\{0\}$, $\pi\colon TM\to M$ the natural  projection, reduces to a linear connection on $M$, even if the Finsler metric is not a Riemannian one. When this happens the Finsler metric is said of {\em Berwald} type.  Since from  a theorem of Szab\'o (cf. \cite[\S 
10.1]{
BaChSh00}), given a Berwald metric $F$ on $M$, there  exists a Riemannian metric  such that its Levi-Civita connection coincides with the Chern connection of $F$, we can state, as already observed in \cite{BorOli10}, that Bishop's proof is also valid in any Berwald space.} 

The approach followed in \cite{BaCaGS11} will be briefly described in the next section. Here, we would like to emphasize that when the Finsler metric comes from a Riemannian one (i.e. $F(v)=\sqrt{g(v,v)}$) the result in \cite{BaCaGS11} improves the one in \cite{Bishop75} with respect to the smoothness assumption on the metric  and on the hypersurface. Indeed from \cite{BaCaGS11},  it is enough that the metric is differentiable with locally Lipschitz differential (and we will assume the same on the semi-Riemannian metric $g$) while in \cite{Bishop75} it assumed that it is $C^4$.  
\bere\label{c2}
In regarding to the smoothness of the hypersurface, we slightly improve the requirement in \cite{BaCaGS11}  that the hypersurface is an embedded submanifold locally described as the regular level of a twice differentiable function having Lipschitz second differential. 
Since at each point of an embedded submanifold $H$ there is a coordinate system adapted to $H$ (see e.g. \cite[Ch.1, Prop. 28]{One83}), it will be clear in the proof of Lemma~\ref{quasiconv}  that it is enough to assume that the hypersurface $H$ is a $C^2$  embedded submanifold of $M$ (of course, also $(M,g)$ must have a differentiable structure of class at least $C^2$).
\ere
\section{Infinitesimal convexity implies local convexity}\label{main}
Let $(M,g)$ be a semi-Riemannian manifold of dimension $m$ and $H$ be an  embedded non-degenerate hypersurface in $M$. Let $\phi$ be a function defined in a neighborhood $\tilde U\subset M$ of  $p_0\in H$ such that $\phi^{-1}(0)= H\cap \tilde U$ and its gradient $\nabla\phi$ at any $p\in  \phi^{-1}(0)$ has the same orientation of the smooth unit normal unit vector field $n$ in 
$H\cap\tilde U$ with respect to which the second fundamental form of $H\cap\tilde U$ is defined. 
Let $H_{\phi}:=\nabla(\de \phi)$ be the Hessian of $\phi$, where $\nabla$  is the Levi-Civita  connection of $(M,g)$. In a coordinate system $(W,(x^1,\ldots,x^m))$  of $M$, using the Einstein summation convention, we get  
\beq(H_{\phi})_p(u,v)=\frac{\partial^2 \phi}{\partial x^i\partial x^j}(p)u^iv^j-\frac{\partial\phi}{\partial x^l}(p)\Gamma^l_{ij}(p)u^iv^j,\label{hess}\eeq
for any $p\in V$ and any $u=(u^1,\ldots, u^m),\ v=(v^1,\ldots,v^m)\in T_pM$, where $\Gamma^l_{ij}$ are the Christoffel symbols of the metric $g$. 
From \eqref{hess}, it's easy to see that if $\gamma\colon (a,b)\to M$ is a geodesic
then
\beq\label{phigamma}
(\phi\circ\gamma)''(s)=(H_{\phi})_{\gamma(s)}\big(\dot\gamma(s),\dot\gamma(s)\big).
\eeq
The following lemma, concerning the equivalence between the second fundamental form of $H\cap\tilde U$ and the Hessian of a function $\phi$ as above, is well known.
We give here a proof for the sake of completeness.
\bl\label{varconv}
For all $V\in T(H\cap\tilde U)$, $H_\phi(V,V)=-|\nabla\phi|\Pi(V,V)$
\el
\begin{proof}
By definition, $H_\phi(V,V)= \tilde V(\tilde V(\phi))-\nabla_{\tilde V}\tilde V(\phi)$,
where $\tilde V$ is any vector field on $M$ extending $V$ and $\nabla\phi$ is the gradient of $\phi$. Using the metric compatibility
of the Levi-Civita connection we get
$\tilde V(\tilde V(\phi))-\nabla_{\tilde V}\tilde V(\phi)=g(\nabla_{\tilde V}\nabla \phi,\tilde V)$.
As $V$ and $\nabla\phi$ are  orthogonal at any point of $H\cap\tilde U$, we get on $H$
\[0=V(g(\nabla\phi,V))=g(\nabla_{\tilde V}\nabla \phi,\tilde V)+g(\nabla\phi,\nabla_{\tilde V}\tilde V),\]
hence $H_\phi(V,V)=-g(\nabla\phi,\nabla_{\tilde V}\tilde V)$ on $H$. As $n=\frac{\nabla\phi}{|\nabla\phi|}|_{H}$ and $\Pi(V,V)=g(n,\nabla_{\tilde V}{\tilde V})$ we  get the thesis.
\qed
\end{proof}
\bere\label{degenerate}
As already observed in \cite[Section 3]{bgs1},  locally defining the hypersurface $H$ as a regular level set of a function allows  to extend the notion of infinitesimal convexity to  hypersurfaces having points where the tensor induced by the ambient metric is degenerate. At these points, the notion of second fundamental form is meaningless. In what follows, unless differently specified, we will assume that $H$ is an embedded hypersurface in $(M,g)$, non necessarily non-degenerate.
\ere
\bd\label{infconvdeg}
Let $(M,g)$ be a semi-Riemannian manifold and $H$ be a $C^2$ embedded hypersurface in $M$. We say that $H$ is infinitesimally convex in a neighborhood $U\subset H$ of a point $p_0\in H$ if there exists a neighborhood $\tilde U$ of $p_0$ in $M$ and a $C^2$ function $\phi\colon \tilde U\to \R$ such that $U=\tilde U\cap H=\phi^{-1}(0)$, $0$ is a regular value of $\phi$, and $H_\phi(V,V)$ is semi-definite (either negative or positive), for all $V\in TU$. 
\ed
\bere\label{phi}
From Lemma~\ref{varconv}, it is clear that this definition is independent of the function $\phi$ if $H$ is a non-degenerate hypersurface. A posteiori, from Theorem~\ref{main}, the same  is true also in the degenerate case. Indeed if there exists a function $\phi$ with respect to which $H$ is infinitesimally convex in $U$, by Theorem~\ref{main}, $H$ is locally convex in $U$ and then it is infinitesimally convex with respect to any other function. Moreover the ``local convex side'' of $H$ (i.e the closure of the connected component of $\tilde U\setminus H$ where the geodesics starting at the points of $U$ with velocity vector tangent to $H$ are locally contained) is $\{x\in \tilde U:\phi(x)\leq 0\}$, if $H_\phi(V,V)\leq 0$, and $\{x\in \tilde U:\phi(x)\geq 0\}$, if 
$H_\phi(V,V)\geq 0$.
\ere
A very natural approach in trying to prove that infinitesimal convexity in a neighborhood of a point $p_0$ of a hypersurface implies local convexity is to evaluate  $\phi$  along any geodesic arc $\gamma$ starting  at $p_0$ with initial velocity vector tangent to $H$. 
The infinitesimal convexity assumption can be used to get a differential  inequality satisfied by the function $\rho=\phi\circ\gamma$,  at least when the image of $\gamma$ is contained on the  side of $H$ which is the candidate convex one (recall Remark~\ref{phi}). The differential inequality and  the initial conditions satisfied by $\rho$ imply that $\rho$ is actually constant and equal to $0$, 
that is the geodesic $\gamma$ is contained in $H$.
This is an intermediate fundamental step to prove local convexity at $p_0$. The remaining part of the  proof is  a quite trivial consequence of the fact that any point in a semi-Riemannian manifold has a convex neighborhood.
\bl\label{quasiconv}
Assume  that $H$ is infinitesimally convex in a neighborhood $U$ of $p_0\in H$. Let $\tilde U$ be a neighborhood of $p_0$ in $M$ and $\phi\colon \tilde U\to\R$ such that $0$ is a regular value of $\phi$ and $\phi^{-1}(0)=H\cap \tilde U=U$, and $H_\phi(V,V)\leq 0$ for each $V\in T U$. Let $\gamma\colon[0,b]\to M$ be a geodesic satisfying the initial conditions $\gamma(0)=p_0$, $\dot\gamma(0)\in T_{p_0}H$ and such that $\gamma([0,b])\subset \{x\in \tilde U\colon \phi(x)\geq 0\}$, then $\gamma([0,b])\subset H$. 
\el
\begin{proof}
Since $H$ is an embedded submanifold, we can assume, without loss of generality, that $\tilde U$ is the domain of a coordinate system, $\varphi\colon \tilde U\to \R^m$, $\varphi=(x^1,\ldots , x^m)$ centered at $p_0$ and that $x^m(p)=0$, for each $p\in H\cap \tilde U$, and the function $\phi$ is  equal to $x^m$. Then the map  $\mathcal P\colon \tilde U\to H$, $\mathcal P(p) =\varphi^{-1}\big((x^1(p),\ldots x^{m-1}(p),0\big)$, is well defined. Let $\gamma_H$ be the projection on $H$ of the geodesic $\gamma$, $\gamma_H(s)=\mathcal P(\gamma(s))=\varphi^{-1}(\gamma^1(s),\ldots,\gamma^{m-1}(s),0)$, where $\gamma^i$, $i=1,\ldots,m$ are the components of the geodesic $\gamma$ in the coordinate system $(\tilde U,\varphi)$ (see Fig. 2). 
\begin{figure}[h]
\includegraphics[scale=1.6]{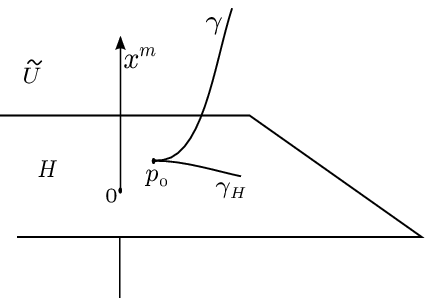}
%
%
\caption{}
\end{figure}
In our notations, the function  $\rho=x^m(\gamma)$ is the $m$-th component $\gamma^m$ of $\gamma$. Thus, from \eqref{phigamma} and \eqref{hess}, the second derivative of $\rho$ is given by
\[\ddot\rho(s)=(H_{x^m})_{\gamma(s)}(\dot\gamma(s),\dot\gamma(s))=-\Gamma^m_{ij} (\gamma(s))\dot\gamma^i(s)\dot\gamma^j(s).\]
Since $H$ is infinitesimally convex in $U$, 
\[(H_{x^m})_{\gamma_H(s)}(\dot\gamma_H(s),\dot\gamma_H(s))=-\Gamma^m_{ij} (\gamma_H(s))\dot\gamma_H^i(s)\dot\gamma_H^j(s)\leq 0.\] 
Thus we can estimate $\ddot\rho$ as follows
\begin{align}
\ddot\rho(s)&\leq -\Gamma^m_{ij}(\gamma(s))\dot\gamma^i(s)\dot\gamma^j(s)+\Gamma^m_{ij} (\gamma_H(s))\dot\gamma_H^i(s)\dot\gamma_H^j(s)\nonumber\\
&=-\Gamma^m_{ij}(\gamma(s))\dot\gamma^i(s)\dot\gamma^j(s)+
\Gamma^m_{ij}(\gamma_H(s))\dot\gamma^i(s)\dot\gamma^j(s)\nonumber\\
&\quad+\Gamma^m_{ij}(\gamma_H(s))(\dot\gamma^i(s)+\dot\gamma_H^i(s))(\dot\gamma_H^j(s)-\dot\gamma^j(s))\nonumber \\
&=\Big(\Gamma^m_{ij}(\gamma_H(s))-\Gamma^m_{ij}(\gamma(s))\Big)\dot\gamma^i(s)\dot\gamma^j(s)\nonumber\\
&\quad +\Gamma^m_{ij}(\gamma_H(s))(\dot\gamma^i(s)+\dot\gamma_H^i(s))(\dot\gamma_H^j(s)-\dot\gamma^j(s))\label{sum}
\end{align}
For each $p\in U$, let us call $\Gamma^m(p)$ the symmetric bilinear operator on $\R^m$ defined by the Christoffel symbols $\Gamma^m_{ij}$ at $p$. The first term in the  summand in \eqref{sum} is bounded above by $||\Gamma^m(\gamma_H(s))-\Gamma^m(\gamma(s))||\,|\dot\gamma(s)|^2$, where $\|\cdot\|$ is the norm of the  bilinear operators  on $\R^m$  and $|\gamma(s)|$ is the Euclidean norm of the   
vector $(\dot\gamma^i(s))_{i}\in \R^m$. Since $\gamma$ is smooth on $[0,b]$ and the Christoffel symbols are Lipschitz on $U$, this last quantity is bounded above by $a|\gamma_H(s)-\gamma(s)|$, where $a$ is a  positive constant depending on $\gamma$ and 
$|\gamma_H(s)-\gamma(s)|$ is the Euclidean norm of $(\gamma_H^i(s)-\gamma^i(s))_{i}\in\R^m$, hence it is equal to $|\gamma^m(s)|=|\rho(s)|=\rho(s)$. By continuity, the second term in  \eqref{sum} is bounded above by  $b |\dot\gamma_H^j(s)-\dot\gamma^j(s)|=b|\dot\gamma^m(s)|=b|\dot\rho(s)|$, where $b$ is a positive constant depending on $\gamma$. Therefore, $\rho$ satisfies the differential inequality 
\beq\label{diffineq}
\ddot\rho\leq c(\rho+|\dot\rho|)
\eeq 
and the initial condition $\rho(0)=\dot\rho(0)=0$ and from \cite[Lemma 3.1]{BaCaGS11}, it must be equal to $0$ on $[0,b]$ (i.e. $x^m\circ\gamma=\phi\circ\gamma\equiv 0$).
\qed
\end{proof}
\bt\label{equiv}
Let $(M,g)$ be a semi-Riemannian manifold and $H$ be an embedded $C^{2}$ hypersurface in $M$. Let $p\in H$ and  $U$ be a neighborhood of $p$ in $H$, then $H$ is infinitesimally convex in $U$  if and only if it is locally convex in $U$.
\et
\begin{proof}
Let $\tilde U$ be an open subset of $M$ and $\phi\colon \tilde U\to \R$ be a function such that $0$ is a regular value, $\phi^{-1}(0)=H\cap \tilde U$ and $U= H\cap \tilde U$. 

Assume that $H$ is locally convex at any point $p\in U$, i.e there exists a neighborhood $O$ of $0\in T_pM$ such that $\exp_p(O\cap T_pH)$ is  contained in the closure one of the connected component of $\tilde U\setminus H$, say
$\{x\in\tilde U\colon \phi(x)\leq 0\}$. Let $v\in O\cap T_pH$ and consider the affinely parametrized geodesic 
$\gamma$ such that $\gamma(0)=p$, $\dot\gamma(0)=v$. Thus there exists $a\in(0,+\infty)$ such that the function $\phi\colon\gamma\colon [-a,1]\to\R $ is well defined and has a maximum point at $0$. Therefore $(H_\phi)_p(v,v)=(\phi\circ\gamma)''(0)\leq 0$.

Assume now that $H$ is infinitesimally convex in $U$, with $H_\phi(V,V)\leq 0$ for all $V\in T U$. Let $p\in U$; we are going to show that  for any open convex neighborhood $B$ of $p$ in $M$ the set $C\colon =B\cap \{x\in \tilde U\colon \phi(x)>0\}$ is also convex, i.e any two points in $C$ are joined by a unique geodesic whose support is contained in $C$.
Let $A$ be the subset of $C\times C$ given by the couple of points that can be connected by a unique geodesic with support in $C$. As $C$ is a connected subset of $M$, it is enough to show that $A$ is non-empty and it is an open and closed subset in $C\times C$. Clearly, each couple $(p_1,p_1)\in C\times C$ can be connected by a constant geodesic, i.e. $A\neq\emptyset$. If $(p_1,p_2)\in A$ and $\alpha$ is the unique geodesic connecting them and whose inner points are in $C$, by smooth dependence of 
geodesics by boundary conditions in a convex 
neighborhood (cf. \cite[Ch. 5, Lemma 9]{One83}), we can consider two small enough neighborhoods $U_1$ and $U_2$ of $p_1$ and, respectively, $p_2$ in $C$, such that the unique geodesic in $B$ connecting $\bar p_1\in U_1$ to $\bar p_2\in U_2$ lies in a small neighborhood of $\alpha$, hence its  points are contained in $C$ and $A$ is open. Now let $(p_1,p_2)\in \overline A\subset C\times C$ and consider a sequence $(p^1_n, p^2_n)_n\subset A$  converging to $(p_1,p_2)$. The sequence of geodesics $\alpha_n$, parametrized on $[0,1]$ and  connecting $p^1_n$ to $p^2_n$ in $B$ converges in the $C^2$ topology to the geodesic $\alpha\colon [0,1]\to B$  connecting $p_1$ to $p_2$. Thus $\alpha$ lies in $\{x\in\tilde U\colon \phi(x)\geq 0\}$. As the points $p_1$ and $p_2$ are in $\{x\in\tilde U\colon \phi(x)>0\}$, from Lemma~\ref{quasiconv}, $\alpha$ cannot be tangent to $H$ at any of its inner points.

Having proved that $C$ is a convex neighborhood in $M$, the rest of the proof follows by contradiction. Indeed, assume that there exists $p\in U$ and a sequence of vectors $(v_n)\subset T_p H$ such that $v_n\to 0$ and $p_n\colon=\exp_p(v_n)\in C$. Let $\gamma_n$ be the geodesic in $C$ connecting the first point of the sequence $p_1$ to $p_n$. As $p_n$ converges to $p$ the sequence $\gamma_n$ converges uniformly (actually in the $C^2$ topology) to the geodesic $\gamma\colon [0,1]\to B$ connecting $p$ to $p_1$ in $B$. By uniform convergence, the image of $\gamma$ must be contained in $\{x\in\tilde U\colon \phi(x)\geq 0\}$. Since $\gamma(0)=p\in H$ and $\dot\gamma(0)=v_1\in T_pH$, from Lemma~\ref{quasiconv}, $\gamma([0,1])\subset H$.  But $\gamma(1)=p_1\not\in H$.
\qed  
\end{proof}
\section{Some remarks and applications}\label{timespace}
\subsection{Convexity with respect to the geodesics having the same causal character}
In a semi-Riemannian manifold, the set of geodesics through a point $p$ can be divided into three disjoint subsets according to the causal character of the initial velocity vector at $p$. It is natural to ask if the equivalence in Theorem~\ref{equiv} holds restricting the geodesics, or equivalently the tangent vectors to $H$ at any $p\in U\subset H$, involved in the definitions of local, or infinitesimal convexity, to one of these subsets. 

To be more precise, let us give the definition of local and infinitesimal convexity, taking into account causality. Let $H$ be a non-degenerate embedded hypersurface of $(M,g)$, $p\in H$ and $W_{p}$ be the set of the timelike (resp. lightlike, spacelike) vectors in $T_{p} M$. Let $U\subset H$ be a neighborhood of $p$ in $H$ and $n$ be a unit normal smooth vector field on $U$.  We say that $H$ is {\em time-} (resp. {\em null-}, {\em space-}) {\em locally convex} at $p_0\in H$ if there exists a neighborhood $O$ of $0\in T_{p_0}H$ such that ${\rm exp}_{p_0}(W_{p_0}\cap O)$ is contained in the closure of the connected component of a tubular neighborhood of $H\cap U$ individuated by $-n$. 
Analogously, we say that a non-degenerate embedded hypersurface $H$ is {\em time-} (resp. {\em null-}, {\em space-}) {\em infinitesimally convex} in a neighborhood $U\subset H$ of $p_0$ if its second fundamental form with respect to $n$
is  positive semidefinite on timelike (resp. null, spacelike) vector field on $TU$, i.e $\Pi_p(v,v)\geq 0$, for each $p\in U$ and for all vectors $v\in W_p\cap T_{p}H$. By Lemma~\ref{varconv}, this last condition is equivalent to $(H_\phi)_p(v,v)\leq 0$, for all  $v\in W_p\cap T_{p}H$, where $\phi$ is a function as in Definition~\ref{infconvdeg}. As in Remark~\ref{degenerate}, defining time- (resp. null-, space-) infinitesimal convexity in this last way allows one to consider also degenerate hypersurface.
\bere\label{nullasboundary}
Since the null cone at a point $p\in M$  is the boundary of both the subsets of spacelike and timelike vectors at $p$,   if $H$ is time- or space-infinitesimally convex at a point $p_0\in H$ then, by  continuity, it is also null-infinitesimally convex. 
\ere
We observe that Lemma~\ref{quasiconv} continues to hold for time- and space-infinitesimal convexities, up to taking a smaller $b>0$ as the upper limit of the interval of definition of the geodesic $\gamma$. Indeed the projection map $\mathcal P$ is smooth and therefore, since $\dot\gamma_H(0)=\dot\gamma(0)$ is timelike (resp. spacelike), the curve $\gamma_H$ remains timelike (resp. spacelike) on a right neighborhood of $0\in\R$, thus inequality \eqref{sum} is true for $s$ in such neighborhood (clearly, for null-infinitesimal convexity, this argument is invalid). 
\bere\label{problema}
On the contrary, under the weaker hypothesis that $H$ is only time- (resp. space-) infinitesimally convex in $U$, the proof of the ``only if'' part of Theorem~\ref{equiv} becomes wrong 
because the set $A$, now defined as the set of the couples of points in $C$ that can be joined by a timelike (resp. spacelike) geodesic contained in $C$ is no longer closed, in fact the limit of a sequence of timelike or spacelike geodesics might be a null geodesic, Remark~\ref{nullasboundary}. Moreover, allowing the limit  be a null geodesic leads to change the definition of $A$ including points that can be connected by a null geodesic but then $A$ becomes non-open. Thus we leave the following as an open problem:
\ere
\begin{problem}
Prove or disprove that for any smooth embedded hypersurface $H$ in a semi-Riemannian manifold (or, at least, in a Lorentzian manifold) time-, null- or space- infinitesimal convexity in a neighborhood $U\subset H$ of a point $p\in H$ implies the same type (time-, null- or space-) of local convexity at $p$.
\end{problem}
\bere\label{nullequiv}
A positive answer to this problem has been given  in \cite[Cor. 3.5 and Rem. 3.8]{CaGeS11} for null/time-infinitesimal and local convexities  
of a hypersurface of the type $H_0\times \R$ in a standard stationary Lorentzian manifold $M_0\times \R$, by reducing the problem to a Finslerian one (cf. \cite{CaJaMa11,CaJaS11} or \cite{Jav}) and using the above mentioned result \cite{BaCaGS11}. 
\ere
\subsection{Geometric convexity}
Let $\Omega$ be an open subset of the semi-Riemannian manifold $(M,g)$ with $C^2$ differentiable boundary $\partial\Omega$. Let $\phi\colon M\to \R$ be a $C^2$ function such that $0$ is a regular value, $\phi^{-1}(0)=\partial \Omega$ and $\phi(x)>0$, for each $x\in\Omega$. 

In what follows, infinitesimal convexity of $\partial \Omega$ will be always considered as defined in terms of the function $\phi$ globally defining $\partial \Omega$ as a regular level set.

Observe that under these assumptions $\partial \Omega$ is a a $C^2$ embededd and oriented hypersuface in $(M,g)$ (the orientation is given by the transversal vector field $\nabla^R \phi$, where $\nabla^R$ is the gradient operator with respect to any auxiliary Riemannian metric on $M$).

From Lemma~\ref{quasiconv}, it immediately follows that if $\partial\Omega$ is infinitesimally convex at any of its points  then the following condition, called (e.g. in  \cite{San01, bgs1, bg, BaCaGS11}) {\em geometric convexity} of $\partial\Omega$ or, equivalently,  of $\Omega$, holds:
\begin{itemize}
\item[$\circ$] for any two points $p,q\in \Omega$ and for any geodesic arc $\gamma\colon[0,1]\to M$ from $p$ to $q$,  if $\gamma([0,1])$ is contained in $\Omega\cup\partial \Omega$ then, actually, it is contained in $\Omega$.
\end{itemize}
(observe that the same also happens for respectively time/space-infinitesimal convexity of $\partial \Omega$ and timelike and spacelike geodesics).

This fact was already observed in \cite[Th. 6]{g} for a $C^3$ domain $\Omega$ in a complete Riemannian manifold. We emphasize that, from a technical point of view,  local convexity at each point of $\partial \Omega$ (which clearly implies geometric convexity), is harder to prove than the other  convexity notions except in the case when {\em strongly} (i.e. \eqref{infconvineq} is satisfied with the strict inequality) infinitesimal convexity holds (cf. \cite[Section 1]{San01}, where strongly infinitesimal convexity is called {\em strictly} infinitesimal convexity); for example, in trying to prove that infinitesimal convexity in a neighborhood of a point $p\in\partial\Omega$ implies local convexity at $p$, we must encompass the possibility of a geodesic $\gamma=\gamma(t)$ oscillating between $\Omega$ and $M\setminus\bar \Omega$, as $t\to t_0$, with $\gamma(t_0)=p$. 
In other words, the fact that $\rho$ is nonnegative is essential in differential inequality \eqref{diffineq}.
Anyway, if geometric convexity holds we immediately get, as in the first part of the proof of Theorem~\ref{equiv}, that $\partial\Omega$ is infinitesimal convex at any of its points. Thus we can state:
\bco\label{allequiv}
For any $C^2$  open subset $\Omega$ of a semi-Riemannian manifold local, infinitesimal and geometric convexity of $\partial \Omega$ are equivalent; moreover time/space-infinitesimal convexity are respectively equivalent to time/space-geometric convexity.
\eco
We recall that the first chain of equivalence holds also for a $C^2$ open subset of a Finsler manifold (see \cite[Cor. 1.2]{BaCaGS11} and recall Remark~\ref{c2}), while both the equivalences between full infinitesimal and geometric convexity and time/space-infinitesimal convexity and time/space-geometric convexity were obtained in \cite[Ths. 4.3, 4.4 and Appendix A]{bgs1} for a $C^3$ open subset of the type
$\Omega=\Omega_0\times \R$ in a standard stationary Lorentzian manifold  $M=M_0\times\R$.
The equivalence in the null case, already obtained in \cite{bgs1} for a static standard region, was proved in \cite[Th. 2.5]{bg}. We observe that this last equivalence also follows by the results in \cite{CaGeS11} as in Remark~\ref{nullequiv}.

Space-, null- and time-geometric convexity, in a Lorentzian setting, were introduced in \cite{bfg} and were used there, together with variational methods, to prove existence and multiplicity results about the number of spacelike and timelike geodesics  connecting a couple of points (for timelike ones, only certain chronologically related points are to be  considered) in an  open subset $\Omega=\Omega_0\times \R$ having space or time-geometrically convex boundary  of a standard static Lorentzian manifold  $M=M_0\times\R$. 
These results were extended to standard stationary Lorentzian manifolds
in \cite{gm}. Remarkably, null- and time-geometric convexity of some open subsets of this type, contained in the outer Schwarzschild, Reissner-Nordstr\"om and Kerr spacetimes, have been proved in \cite[\S 7]{mas} (see also \cite{bfg2,FlSa,CaGeS11}). We stress that these results were obtained by showing that the boundaries of such open subsets are strongly infinitesimally convex and indeed the advantage of the notion of infinitesimal convexity with respect to the other ones relies, of course, in its direct computability.
\subsection{Some applications to geodesic connectedness}
It is worth to observe that geometric convexity extends the classical notion of convexity of a subset of $\R^m$. For example, assume that 
$(M, g)$ is a smooth complete Riemannian manifold and $\Omega$ is a smooth, connected open subset of $M$ having geometrically convex boundary $\partial\Omega$;
then there exists a  (non necessarily unique) geodesic connecting $p$ to $q$,  contained in $\Omega$ and having length equal to the distance in $\Omega$ between $p$ and $q$. A way to prove that  is to apply both a minimization and a penalization argument to the energy functional of the Riemannian manifold. We refer to  \cite[Cor. 4.4.7]{mas} for a proof, using these variational methods, of the existence of a geodesic connecting $p$ to $q$; the minimizing property of such a geodesic can be proved as in \cite[Remark at p.448]{g1} while
the case when $\partial \Omega$ is not differentiable and/or  $(M,g)$ is not complete has been studied in \cite{bgs1}.

This result also holds in a forward or backward complete Finsler manifold \cite[Th. 1.3]{BaCaGS11}. A more general result is obtained replacing the completeness of
the Finsler manifold $(M,F)$ with the assumption that the closure in $\Omega$ of any ball, with respect to the symmetric distance associated to the the  pseudo-distance induced by $F$ on $\Omega$, is compact \cite[Th. 1.3]{BaCaGS11} and \cite[Rem. 4.2]{CaGeS11}.

Another way to prove such results is by using a shortening argument to the length functional.
In an open subset of a  Lorentzian manifold satisfying good causality properties,  similar and in some aspects dual (cf. \cite[p. 409]{One83}) arguments are applicable.  We refer to \cite[Ch. 14]{One83} for more details on the the definitions and notations about causality that we are going to use.

Let $(M,g)$ be a time-oriented Lorentzian manifold and $\gamma$ be a causal curve on $M$ (i.e., assuming for simplicity that $\gamma$ is piecewise smooth, $g(\dot\gamma,\dot\gamma)\leq 0)$). The Lorentzian length of $\gamma$ is defined as 
$L(\gamma)=\int_{\gamma}\sqrt{-g(\dot\gamma,\dot\gamma)}$. Let $\Omega$ be a subset of $M$ and $p$, $q$ two points in $\Omega$. We will denote by $\mathcal{C}^{\Omega}_{pq}$ the set of the future pointing causal curves $\gamma\colon[0,1]\to \Omega$ such that $\gamma(0)=p$ and $\gamma(1)=q$.  We say that
{\em $q$ is causally related to $p$ in $\Omega$} and we  write $p<^{\Omega}q$, if $\mathcal{C}^{\Omega}_{pq}\neq \emptyset$. We define the {\em time separation in $\Omega$}, $\tau^{\Omega}(p,q)\in[0,+\infty]$, as $\tau^{\Omega}(p,q)=\sup_{\gamma\in \mathcal{C}^{\Omega}_{pq}} L(\gamma)$, if $\mathcal{C}^{\Omega}_{pq}\neq\emptyset$ otherwise
$\tau^{\Omega}(p,q)=0$.
\bpr\label{causalcompact}
Let $(M,g)$ be a time-oriented Lorentzian manifold and $\Omega$ be a $C^2$,
open subset such that $\bar\Omega$ is compact and strongly causal. Assume that $\partial\Omega$ is convex (recall the first part of Corollary~\ref{allequiv}), then for any $p,q\in \Omega$ with $p<^{\Omega}q$   there is a future pointing causal geodesic  $\gamma\in \mathcal{C}^{\Omega}_{pq}$ such that $L(\gamma)=\tau^{\Omega}(p,q)$ (hence $\tau^{\Omega}(p,q)$ is finite).
\epr
\begin{proof}
Let $\{\gamma_n\}$ be a sequence of curves, $\gamma_n:[0,1]\to\Omega$, such that $\gamma_n(0)=p$, $\gamma_n(1)=q$ and $L(\gamma_n)\to \tau^{\Omega}(p,q)$. From  \cite[Ch. 14, Prop. 8]{One83} there is a  limit sequence $p_0=p<^M p_1<^M\ldots<p_j<^M\ldots$ (see \cite[Ch. 14, Def. 7]{One83}) for $\{\gamma_n\}$. Arguing as in the proof of \cite[Ch.  14, Lemma 14]{One83} the limit sequence must be finite, and the last point $p_m$ has to  be equal to $q$.
From the definition of a limit sequence, there exist a subsequence  $\{\gamma_k\}$ 
and sequences $\{s_{kj}\}$,  $j=0,\ldots,m$, such that $\gamma_k(s_{kj})\to p_j$ and then $p_j\in\bar\Omega$; moreover there are  convex neighborhoods $\{C_j\}_{j=0,\ldots,m-1}$ associated to the limit sequence $\{p_j\}_{j=0,\ldots,m}$ such that $p_j,p_{j+1}\in C_j$, for each $j=0,\ldots,m$.   From the second part of the proof of Theorem~\ref{equiv}, we know that the sets $C_j\cap\Omega$ are  convex neighborhoods. Consider the  future pointing causal broken geodesic $\lambda$ from $p$ to $q$, with vertices $p_j$, having one segment $\lambda_j$ in each convex set $C_j$ (such broken geodesic is a {\em quasi-limit} of $\{\gamma_n\}$, \cite[p. 406]{One83}). 
Consider also the future pointing causal geodesics $\lambda_{kj}$ connecting $\gamma_k(s_{kj})$ with $\gamma_k\big(s_{k(j+1)}\big)$ in\footnote{Observe that such geodesics must be causal and future pointing as the points  $\gamma_k(s_{kj})$ and $\gamma_k\big(s_{k(j+1)}\big)$ are causally related, see \cite[Ch. 14, Lemma 2]{One83}.} $C_j$.
By the smooth dependence of geodesics by boundary conditions in a convex neighborhood (\cite[Ch.5, Lemma 9]{One83}) the geodesics $\lambda_{kj}$ converge in the $C^2$ topology to $\lambda_j$ and $L(\lambda_{kj})\to L(\lambda_j)$.  
Hence $\lambda$ is contained in $\bar\Omega$ and, by the length maximizing property of each causal geodesic segment in $C_j$ (see \cite[Ch. 5, Prop 34]{One83}), we get  $\tau^{\Omega}(p,q)=\lim_kL(\gamma_k)\leq \lim_k L(\lambda_k)=L(\lambda)$, where $\lambda_k$ is the broken geodesics in $\Omega$ with vertices $\gamma_k(s_{kj})$ and segments $\lambda_{kj}$. Since $p$  belongs to $\Omega$, the portion of $\lambda$ given by the first segment $\lambda_0$ plus  a small part   $\lambda^s_{1}\subset\lambda_1$ starting at $p_1$ and contained in $C_0$ cannot be entirely contained in $\partial\Omega$ and then, by Lemma~\ref{quasiconv}, it cannot be a pregeodesic. Therefore there exists a causal future pointing geodesic $\alpha_0$ having length strictly
greater than $L(\lambda_0)+L(\lambda^s_{1})$, connecting $p$ to the end point $p_1^s$ of $\lambda^s_1$.
By choosing a sequence of points $\{p^s_{h1}\}\subset C_0\cap\Omega$, such that $p^s_{h1}\to p^s_{1}$, the future pointing causal geodesics defined by the points $p_0$ and $p^s_{h1}$ and having support in $C_0\cap\Omega$ converge, in the $C^2$ topology, to $\alpha_0$. Hence by Lemma~\ref{quasiconv}, $\alpha_0$  does not intersect $\partial\Omega$ in any of its inner points.  As the end point $p_1^s$ belongs to $C_0\cap C_1$, in an analogous way, replacing $C_0$ with $C_1$, we can construct a  broken geodesic segment $\alpha_2$ between $p$ and $p_2$, longer than $\alpha_0$ plus the segment of $\lambda$ from $p_1^s$ and $p_2$, which does not intersect $\partial\Omega$ in any of its inner point. Thus in $m$ steps, we  construct a future pointing causal broken  geodesic $\alpha$ connecting $p$ to $q$, entirely contained in $\Omega$ and such that $\tau^{\Omega}(p,q)\leq L(\lambda)<L(\alpha)$. This contradiction comes from the fact that we have assumed that $\lambda$ is not a geodesic. 
Using again Lemma~\ref{quasiconv}, we conclude that $\lambda$ must be contained in $\Omega$ and $\tau^{\Omega}(p,q)= L(\lambda)$.
\qed
\end{proof}
Since a globally hyperbolic Lorentzian manifold is strongly causal and the intersection of the causal future $J^+(p)$ of $p$ with the causal past $J^-(q)$ of $q$ is compact, as in  Proposition~\ref{causalcompact} we get the following.
\bco\label{globconv}
Let $(M,g)$ be a globally hyperbolic Lorentzian manifold and $\Omega$ be a 
$C^2$ open subset. Assume that $\partial\Omega$ is convex, then for any $p,q\in \Omega$ with $p<^{\Omega}q$   there is a future pointing causal geodesic  $\gamma\in \mathcal{C}^{\Omega}_{pq}$ such that $L(\gamma)=\tau^{\Omega}(p,q)$.
\eco
We don't know if we can replace the assumption of   (infinitesimal or, equivalently, local) convexity of $\partial\Omega$ in Proposition~\ref{causalcompact} with  time-infinitesimal convexity (or even time and null local convexity). In fact, in that case, we cannot state that the sets $C_i\cap \Omega$ associated to the limit sequence in the proof of Proposition~\ref{causalcompact} are convex
neighborhoods (recall Remark~\ref{problema}).
\begin{problem}\label{problema2}
Does the conclusion of Proposition~\ref{causalcompact} (or of Corollary~\ref{globconv}) hold assuming that $\partial\Omega$ is time-infinitesimally convex? 
\end{problem}
In some globally hyperbolic stationary Lorentzian manifold we can obtain the full geodesic connectedness of an open  connected subset having convex boundary. 
Indeed  the following proposition holds.
\bpr
Let $(M,g)$ be a  Lorentzian manifold endowed with a  complete timelike Killing vector field $K$ and admitting a smooth, spacelike, complete Cauchy hypersurface. Let $\Omega$ be a $C^2$ open connected subset of $M$ having convex boundary and which is invariant for the flow of $K$. Then for each $p$ and $q$ in $\Omega$ there is a geodesic connecting $p$ to $q$ whose support is in $\Omega$. Moreover, if $\Omega$ is also non contractible, a sequence $\{\gamma_n\}$ of spacelike  geodesics in $\Omega$ connecting $p$ to $q$ exists such that $g(\dot\gamma_n,\dot\gamma_n)\to+\infty$ and, if $l\colon \R\to M$ is a flow line of $K$ with  $l(\R)\subset\Omega$, then the number of timelike future pointing geodesics in $\Omega$ connecting each $p\in\Omega$  to $l(s)$ diverges as $s\to +\infty$.   
\epr
The above proposition extends the results in \cite{piccio}, where $K$ is also assumed to be static (i.e. its orthogonal distribution is integrable).
It is based on the fact that the existence of a smooth, spacelike, Cauchy hypersurface and the completeness of $K$ imply that $(M,g)$ is isometric by the map  $\psi\colon S\times\R\to M$, where $\psi$ is the flow   of $K$ and $S$ the Cauchy hypersurface,
to a standard stationary spacetime  $S\times \R$ (cf. \cite[Th. 2.3]{CFSadvances}). As $\Omega$ is invariant for the flow of $K$, then  $\Omega=\Omega_0\times \R$, with $\Omega_0=S\cap \Omega$ and since $S$ is complete, $\Omega_0\cup\partial \Omega_0$ is complete.  Then, as in the case without boundary in \cite{CFSadvances}, one can prove that the projections on $\Omega_0$ of the curves in a sublevel of the restriction of the energy functional of the standard stationary region $\Omega_0\times\R$ to the manifold of the $H^1$ curves $\gamma$ connecting $p$ to $q$ and such that $g(\dot\gamma, K)$ is constant a.e., are contained in a compact subset of $\Omega_0\cup \partial\Omega_0$. Then, using a penalization argument as in \cite[Section 4.2]{mas},  the result follows as in \cite[Sections 4.4 and 4.5]{mas}.

We point out that, due to the particular variational setting,   last part of the statement of the above Proposition  follows as well  assuming that $\partial\Omega$ is only time-infinitesimally convex (recall Problem~\ref{problema2}). For a standard stationary region $\Omega_0\times\R$, a result of this type, without assuming global hyperbolicity of the ambient spacetime  $(M,g)$, has been obtained in \cite[Th. 4.6]{CaGeS11}.
\subsubsection*{Acknowledgements}
I would like to thank M. S\'anchez for suggesting to investigate the topic of this note and for several useful comments. Moreover, I thank the local organizing committee of the ``VI International Meeting on Lorentzian Geometry, Granada 2011'' for the financial support and the hospitality during the workshop.

\end{document}